\documentclass{amsart}

\usepackage{amssymb}
\usepackage{amsmath}
\usepackage{amsfonts}
\usepackage{natbib}

\usepackage{picins}

\usepackage[all]{xy}
\xyoption{arc} \xyoption{poly}

\newdir{ >}{{}*!/-10pt/@{>}} 

\newtheorem*{lemma*}{Lemma}

\newcommand{\fr}[1]{\ensuremath{\mathfrak{#1}}} 

\begin{document}

\renewcommand{\thefootnote}{\fnsymbol{footnote}}

\begin{center} \Large\textbf{The finiteness theorem for invariants of a finite group}\normalsize\vspace{2ex}

by

Emmy Noether in Erlangen (trans. Colin McLarty)
\end{center}

\vspace{4ex}

The following is an entirely elementary finiteness proof---using
only the theory of symmetric functions---for the invariants of a
\emph{finite} group, which at once supplies \emph{an actual statement of a complete
system of invariants} while the usual proof using the Hilbert basis
theorem (\textbf{Ann}.\ 36) is only an existence proof.\footnote{)
See for example Weber, \emph{Lehrbuch der Algebra} (2.\ Aufl.) 2.\ Band,
\S 57.})

Let the finite group \fr{H} consist of $h$ linear transformations
(with non-vanishing determinant) $A_1\cdots A_h$ where $A_k$
represents the linear transformation \[x_1^{(k)}= \sum_{\nu=1}^n
a_{1\, \nu}^{(k)}x_{\nu},\cdots,x_n^{(k)}= \sum_{\nu=1}^n a_{n\,
\nu}^{(k)}x_{\nu}\] or more briefly: $(x^{(k)})=A_k (x)$.  So the
group \fr{H} takes the series $(x)$ with elements $x_1\cdots x_n$
over into the series $(x^{(k)})$ with elements $x_1^{(k)}\cdots
x_n^{(k)}$.  Since the identity is among the $A_1\cdots A_h$ so
$(x)$ is among the series $(x^{(k)})$. An integral rational
(absolute) invariant of the group is understood to be an integral
rational function of $x_1\cdots x_n$ which is unaltered by
application of the $A_1\cdots A_h$; so that when $f(x)$ is such an
invariant:\setcounter {equation}{0}
\begin{equation}
    f(x)=f(x^{(1)})=\cdots = f(x^{(k)}) =
\frac{1}{h}\sum_{k=1}^{h} f(x^{(k)})
\end{equation}

1.\ Formula (1) says $f(x)$ is an integral, rational,
\emph{symmetric} function of the series of quantities $(x^{(1)})\cdots
(x^{(h)})$ and in fact, because each summand $f(x^{(k)})$ contains
only the one series of quantities $(x^{(k)})$ it is the simplest case,
customarily called \emph{uniform}.  By the well known theorem
on symmetric functions of \emph{series} of quantities\footnote{) See the
footnote to section 2.}) then $f(x)$ has an integral rational
presentation by the elementary symmetric functions of these series,
i.e.\ by the coefficients $G_{\alpha\, \alpha_1 \,\cdots \alpha_n}(x)$ of
the ``Galois resolvent'':
\[ \Phi (z,u) = \prod_{k=1}^h (z+u_1 x_1^{(k)}+\dots + u_n
x_n^{(k)})\]
\[= z^h+ \sum G_{\alpha\, \alpha_1 \,\cdots \alpha_n}(x)z^{\alpha}u_1^{\alpha_1}\cdots u_n^{\alpha_n}
        \binom {(\alpha+\alpha_1 +\dots+ \alpha_n)=h} {\alpha \neq h}\]
where the  $G_{\alpha\, \alpha_1 \, \alpha_n}(x)$ are invariants of
degree $\alpha_1 +\dots+ \alpha_n$ in $x$. So is proven:

\emph{The coefficients $G_{\alpha\, \alpha_1 \, \alpha_n}(x)$ of the
Galois resolvent make up a complete system of invariants of the
group so every invariant has an integral rational presentation
in terms of these finitely many invariants.}

2. One can also use the following even more elementary considerations,\footnote{) This results from a proof of the theorem on symmetric functions of series of quantities in the above mentioned \emph{uniform} case.}) based on (1), which at once yield a second complete system.  Let
\[ f(x) = a + b x_1^{\mu_1}\cdots x_n^{\mu_n}+ \cdots + c x_1^{\nu_1}\cdots x_n^{\nu_n}\]
where $a,b,\cdots, c$ refer to constants.  Then according to (1):
\[ h\cdot f(x) = h\cdot a + b\cdot \sum_{k=1}^{h} x_1^{(k)^{\mu_1}}\cdots x_n^{(k)^{\mu_n}}+ \cdots + c\cdot \sum_{k=1}^{h} x_1^{(k)^{\nu_1}}\cdots x_n^{(k)^{\nu_n}}\]
\[ = h\cdot a + b\cdot J_{\mu_1\cdots \mu_n}+ \cdots + c\cdot  J_{\nu_1\cdots \nu_n}\]
So every invariant is a integral linear combination of the special:
\[J_{\mu_1\cdots \mu_n} = \sum_{k=1}^{h} x_1^{(k)^{\mu_1}}\cdots x_n^{(k)^{\mu_n}}\]and it suffices to give the proof for these.  But, apart from a numerical factor, $J_{\mu_1\cdots \mu_n}$ is the coefficient of $u_1^{\mu_1}\cdots u_n^{\mu_n}$ in the expression \[ S_{\mu} = \sum_{k=1}^{h} (u_1x_1^{(k)}+ \cdots + u_n x_n^{(k)})^{\mu}\] where $\mu = \mu_1+\cdots + \mu_n$, which represents the $\mu$-th power sum of the $h$ linear forms \[\xi_1 = u_1x_1^{(1)}+ \cdots + u_n x_n^{(1)},\cdots, \xi_h = u_1x_1^{(h)}+ \cdots + u_n x_n^{(h)}\] Now it is known that the infinitely many power sums $S_{\mu}$ are all integral rational combinations of \[S_1,\cdots,S_h\] whose coefficients are given by the invariants
\[J_{\mu_1\cdots \mu_n}\hspace{10ex} \mu_1+\cdots +\mu_n \leqq h\] so a second complete system is given:

\emph{A full system of invariants for the group is given by the invariants $J_{\mu_1\cdots \mu_n}$ where $\mu_1+\cdots +\mu_n \leqq h$ and $h$ is the order of the group.}

The connection with 1.~is shown by the remark that when the power sums $S_1\cdots S_h$ are given in terms of the elementary symmetric functions of $\xi_1\cdots \xi_h$ this leads to to the given complete system $G_{\alpha\, \alpha_1 \, \alpha_n}(x)$.  Both results show \emph{all invariants have an integral rational expression in terms of those whose degree in $x$ does not exceed the order $h$ of the group}.

3. These results have consequences for rational representations.  Every rational absolute invariant can be represented---as shown when both the numerator and denominator are expanded by all conjugates of the denominator under the group action---as a quotient of two integral rational absolute invariants which need not be free of common factors.  It follows that \emph{every rational absolute invariant can be expressed rationally by the coefficients  $G_{\alpha\, \alpha_1 \, \alpha_n}(x)$ of the Galois resolvent $\Phi (z,u)$}.  This theorem can easily be proved another way without going through the finiteness theorem; as already in Weber II \S 58.\footnote{)That proof is incorrect.  It only shows the function $\Psi(t)$ in formula (7) has invariant coefficients without showing these invariants can be rationally expressed by coefficients of $\Phi(t)$.  This gap is avoided by using a well known differentiation process instead of the Lagrange interpolation formula to represent the $x_i^{(k)}$ by the $\xi_i$.  Differentiating the identity $\Phi(-\xi_k,u)=0$ at each $u$ gives a relation:
\[ \left[ \frac{\partial \Phi}{\partial u_i} - x_i^{(k)}\cdot \frac{\partial \Phi}{\partial z} \right]_{z=-\xi_k} = 0\] And so, instead of formulas (7) and (8) of Weber for each rational function $\omega(x)$, the following occurs:
\[ \omega (x_1^{(k)}\cdots x_n^{(k)}) = \omega\! \left( \frac{\frac{\partial \Phi}{\partial u_1}}{\frac{\partial \Phi}{\partial z}}\cdots  \frac{\frac{\partial \Phi}{\partial u_n}}{\frac{\partial \Phi}{\partial z}} \right)_{z=-\xi_k}\] which involves only coefficients of $\Phi(z,u)$ and whose summation over all $k$ gives the desired representations of invariants $\omega (x)$.})

4.  In closing let it be said that the results proved here are implicit in my paper ``K\"orper und Systeme rationaler Funktionen"\footnote{) Math. Ann. \textbf{76}, p.~161 (1915).}); specifically the complete system of 1.~in Theorems VI and VII, and a proof for rational invariants independently of the finiteness theorem in Theorem III.\footnote{)Theorems VIII and IX contain a finiteness proof for \emph{relative} invariants which rests on different grounds than the usual but which is still only an existence proof; this is because relative invariants do not form a field.})  The proof methods and results are much simplified in the special case here compared to the general theory.  I came to apply the general investigation to invariants of finite groups through conversations with Herr E.~Fischer who is also the source for the substance of the proof in 2---the complete system there is not rationally definable in the general case---and for the note to 3.\vspace{3ex}

Erlangen, May 1915

\end{document}